\documentstyle[12pt]{article}
\newtheorem{cor}{Corollary}[section]
\newtheorem{prop}{Proposition}[section]
\newtheorem{theo}{Theorem}[section]
\newtheorem{lem}{Lemma}[section]
\newcounter{bean}
\addtocounter{section}{-1}
\begin{document}
\addtolength{\baselineskip}{6pt}
\title{Recovering of curves with involution by extended Prym data}
\author{Vassil Kanev\footnotemark[1]}
\date{  }
\maketitle
\section*{Introduction}
\footnotetext{Supported in part by the Bulgarian foundation
"Scientific research" and by NSF under the US-Bulgarian project
"Algebra and algebraic geometry".}
The classical Torelli theorem states that every smooth, projective
algebraic curve $X$ is determined uniquely up to isomorphism by its
principally polarized Jacobian $(J(X),\Theta )$. In this paper we
consider curves $\tilde{C}$ with involution without fixed points
$\sigma :\tilde{C}\longrightarrow \tilde{C}$. We let
$C=\tilde{C}/\sigma $ and denote by $\pi :\tilde{C}\longrightarrow
C$ the factor map. One associates the principally polarized Prym
variety $(P,\Xi )$ where $P=P(\tilde{C},\sigma )=(1-\sigma
)J(\tilde{C})$ and $\Theta \mid _P$ is algebraically equivalent to
$2\Xi$. The natural question,
whether the Prym variety determines uniquely the pair
$(\tilde{C},\sigma )$ up to isomorphism, has negative answer in
general if the genus $g$ of $C$ is $\leq 6$ as well as in every genus
 $\geq 7$ for some special loci of curves, e.g for hyperelliptic $C$.
The problem of spotting the pairs $(\tilde{C},\sigma )$ which are not
determined uniquely by the Prym variety is still open. Some
partial results have been obtained in
\cite{f-s},\cite{kan},\cite{don},\cite{don1},\cite{deb1},\cite{deb2},
\cite{nar},\cite{verra}.

    In this paper we propose an extension of the Prym data $(P,\Xi )$
 and prove that it determines uniquely up to isomorphism any pair
$(\tilde{C},\sigma )$ for $g\geq 2$. Our extension originates from the
 following observation. Consider the case of $\tilde{C}$ of genus 1.
Here the involution is a translation $t_{\mu }$ by some point $\mu $
of order 2 in $J(\tilde{C})$. Notice that $\{ 0,\mu \} =Ker(Nm_{\pi }:
J(\tilde{C})\longrightarrow J(C))$. The Prym variety is equal to 0.
There is a classically known data which determines uniquely the pair
$(\tilde{C},t_{\mu })$ up to isomorphism (see e.g. \cite{mum2},
\cite{clem}). Namely, one can always represent $J(\tilde{C})$ as
${\bf C}/{\bf Z}\tau +{\bf Z}$, where $\tau $ belongs to the
upper-half plane ${\cal H}$, so that
$$
\mu =\frac{1}{2}\tau +\frac{1}{2}(mod {\bf Z}\tau +{\bf Z}).
$$
So, the moduli space of pairs $(\tilde{C},\sigma )$ is isomorphic to
$\Gamma _{1,2}\backslash {\cal H}$ where $\Gamma _{1,2}\subset
PGL(2,{\bf Z})$ is the subgroup
\[
\Gamma _{1,2}=\{ \left( \begin{array}{cc}a&b\\c&d \end{array} \right)
 \mid ad-bc=1,ab\equiv 0(mod 2),cd\equiv 0(mod 2)
\} \]
One considers the three even theta functions with characteristics :
$\theta _{00}(z,\tau ),\theta _{01}(z,\tau )$ and $\theta
_{10}(z,\tau )$. Just one of them vanishes on $\mu $, namely $\theta
_{00}(z,\tau )$. Now, by the transformation law for theta functions
one checks that for
\begin{equation}\label{ii.1}
\lambda (\tau )=-\theta _{01}(0,\tau )^4/\theta _{10}(0,\tau )^4
\end{equation}
the set $\{ \lambda (\tau ),1/\lambda (\tau )\} $, or equivalently
the function $k(\tau )=\lambda (\tau )+1/\lambda (\tau )$, remain
invariant with respect to the action of $\Gamma _{1,2}$. Moreover the
map
$$
k:\Gamma _{1,2}\backslash {\cal H}\longrightarrow {\bf C}^*-\{ 0,2\}
$$
is an analytic isomorphism (see Section~(\ref{s2})).

    Generalizing to genus greater then 1, first we have that
$Ker(Nm_{\pi }:J(\tilde{C})\longrightarrow J(C))=P\cup P_{\_ }$ where
$P_{\_ }$ is a translation of $P$ by a point of order 2. One considers
 the symmetric theta divisors $\Theta $ of $J(\tilde{C})$ which have
the property that for every $\rho \in P_2$ either $\rho \not \in
\Theta $ or $\rho \in \Theta $ and $mult_{\rho }\Theta $ is even. It
turns out that there are three $P_2$-orbits of symmetric theta
divisors with this property. The divisors of one of the orbits contain
 $P_{\_ }$ ; these are exactly the theta divisors which appear in
Wirtinger's theorem \cite{fay},\cite{mum} and satisfy $\Theta .P=2\Xi
$ for symmetric theta divisors $\Xi $ of $P$. None of the divisors of
the other two orbits contains $P_{\_ }$. Restricting the latter to $P$
 we obtain two $P_2$-orbits of divisors in the linear system
$\mid 2\Xi \mid $ which we denote by $O_1,O_2$.

    The extended Prym data consists of $(P,\Xi )$ together with the
two $P_2$-orbits $O_1,O_2\subset \mid 2\Xi \mid $. Our result is that
for $g\geq 2$ it determines uniquely up to isomorphism any pair
$(\tilde{C},\sigma )$. The proof is analogous to Andreotti's proof
of Torelli's theorem and uses the Gauss map for the divisors of $O_i$.
 Special treatment is required if $C$ is hyperelliptic, or
bi-elliptic,
or g=3. In fact the bi-elliptic case has been already considered
earlier by Naranjo in \cite{nar} who proved that for $g\geq 10$ the
pair $(\tilde{C},\sigma )$ can be recovered by the ordinary Prym data
$(P,\Xi )$. His arguments however do not work for $g=4$ or $5$. If
$\eta $ is the point of order 2 in $Pic^0(C)$ which determines the
covering $\pi:\tilde{C}\longrightarrow C$ then the divisors of $O_i$
are equal to translations of connected components of the set
$$
Z=\{ L\in Pic^{2g-2}(\tilde{C})\mid Nm_{\pi }(L)\simeq K_C\otimes
\eta ,\; h^0(\tilde{C},L)\geq 1\}
$$
It is interesting that line bundles $L$ of this type appear also in
the study of rank 2 vector bundles on $C$ with canonical determinant
\cite{be} as well as in representing $(\tilde{C},\sigma )$ as the
spectral curves associated to $sp(2n)$ -matrices with parameter.

   {\bf Acknowledgement.} Part of this work was done while the author
was a visitor in the University of Michigan and the University of
Utah. The hospitality of these institutions is gratefully
acknowledged.
\begin{center}{\bf Contents}\end{center}
0.Notation and preliminaries.  1.Double unramified coverings of
elliptic curves.  2.Extended Prym data.  3.The semicanonical map and
the Gauss map.  4.The hyperelliptic case, $g\geq 2$.  5.The
bi-elliptic case, $g\geq 4$.  6.  The case $g=3$.

\section{Notation and preliminaries}\label{s1}
We denote by $\equiv$ the linear equivalence of divisors. Let $X$ be
an algebraic, smooth, irreducible curve. We denote by $J_d(X)$ the
divisor classes of degree $d$ modulo linear equivalence and by
$Pic^d(X)$ the isomorphism classes of invertible sheafs of degree $d$
on $X$. Abusing the notation we shall write by the same letter an
element in $J_d(X)$ and the corresponding element in $Pic^d(X)$. If
$D$ is a divisor of degree $d$ and $\xi  = cl(D)$ its class in
$J_d(X)$ we write
$$
h^0(C,D) = h^0(C,\xi  ) = dim\mid \xi  \mid + 1
$$
If $L$ is an invertible sheaf of $X$, then $\mid L \mid $
is the linear system of divisors of sections of $H^0(X,L)$. If $\mid
L \mid $ is without base points we denote by $\varphi  _L = \varphi
_{\mid L
\mid } $ the map $X \longrightarrow \mid L \mid ^*$ defined by
$\varphi _L(x) = \mid L(-x) \mid +x$.

Let $A$ be  a principally polarized abelian variety of dimension $g$
isomorphic to ${\bf C}^g/\Lambda _{\tau }$, where $\Lambda _{\tau } =
{\bf Z}^g\tau  + {\bf Z}^g$ with $\tau \in {\cal H}_g $, where
${\cal H}_g$ is the Siegel upper-half space. Any point $e \in
{\cal C}^g$ has two characteristics $\epsilon ,\delta \in {\bf R}^g$
such that $e =\epsilon \tau + \delta $. We shall write sometimes $e =
\left[ \begin{array}{c} \epsilon \\ \delta \end{array} \right] $. Let
$x \in
 A$ and let $x = (x'\tau +x'')(mod \Lambda _{\tau})$ with $x',x''\in
 {\bf R}^g$. If no confusion arises we shall refer to $x',x''$ as the
characteristics of $x$, keeping in mind that $x',x''$ are determined
modulo ${\bf Z}^g$. We shall denote by $A_2$ the points of order 2 in
$A$. Let $\lambda =\lambda '\tau +\lambda '' ,
 \mu = \mu '\tau + \mu '' \in
A_2 $. The Weyl pairing $e_2 : A_2 \times A_2 \longrightarrow
{\bf Z}$ is defined by
$$
e_2(\lambda ,\mu) = 4(\lambda '^t\mu '' - \mu '^t\lambda ''(mod 2)
$$
Let $\Theta $ be a symmetric theta divisor. One defines a
quadratic form $q_{\Theta } : A_2 \longrightarrow {\bf Z}_2 $
associated with $\Theta $ by
\begin{equation}\label{e5.1}
q_{\Theta }(\lambda ) = mult_0(\Theta +\lambda )(mod 2)
\end{equation}
Consider the theta function with characteristics $\lambda ',\lambda
''\in \frac{1}{2} {\bf Z}^g$
$$
\theta \left[ \begin{array}{c} \lambda '\\ \lambda'' \end{array}
\right] (z,\tau )
= \sum_{n\in {\bf Z}^g}\exp {(\pi i(n+\lambda ')\tau ^t(n+\lambda ') +
2\pi i(n+\lambda '')^t(z+\lambda ''))}
$$
Then
$$
\theta \left[ \begin{array}{c} \lambda '\\ \lambda'' \end{array}
\right] (-z,\tau ) =
\theta \left[ \begin{array}{c} -\lambda '\\ -\lambda'' \end{array}
\right] (z,\tau ) =
(-1)^{4\lambda '^t\lambda ''}
\theta \left[ \begin{array}{c} \lambda '\\ \lambda'' \end{array}
\right] (z,\tau )
$$
Thus if $\Theta $ is the divisor of the theta function
$\theta \left[ \begin{array}{c} 0\\ 0 \end{array} \right] (z,\tau )$
then for any $\lambda =\lambda '\tau +\lambda '' \in A_2$ one has
\begin{equation}\label{e6.2}
q_{\Theta }(\lambda ) = 4\lambda '^t\lambda ''(mod 2)
\end{equation}
For any symmetric theta divisor $\Theta $ of $A$ the bilinear form
associated with $q_{\Theta }$ is $e_2$, i.e.
$$
\Delta _{\lambda }\Delta_{\mu}q_{\Theta }(\xi )
= q_{\Theta }(\xi +\lambda +\mu)-q_{\Theta }(\xi +\lambda )
-q_{\Theta }(\xi +\mu )+q_{\Theta }(\xi )
$$
is independent of $\xi $ and equals $e_2(\lambda ,\mu )$. In
particular one has the following formula
\begin{equation}\label{e6.1}
q_{\Theta }(\lambda +\alpha ) = q_{\Theta }(\lambda ) + e_2(\lambda
,\alpha ) + q_{\Theta }(\alpha ) - q_{\Theta }(0)
\end{equation}
The map $\Theta \longmapsto q_{\Theta }$ gives a bijective
correspondence between the set of symmetric theta divisors of $A$ and
the set of quadratic forms on $A_2$ whose associated bilinear form
equals $e_2$ and which vanish on $2^{g-1}(2^g+1)$ points of $A_2$.

Let $\tilde{C}$ be an algebraic, projective, smooth, irreducible
curve. Let $\sigma : \tilde{C} \longrightarrow \tilde{C}$ be an
involution without fixed points, let $C= \tilde{C}/\sigma$ and let
$\pi : \tilde{C} \longrightarrow  C$ be the factor map. Let $g = g(C),
 \tilde{g} = g(\tilde{C}) = 2g-1 $. We suppose that $g\geq2$. Let
$\tilde{J} = J(\tilde{C}) , J= J(C) $ be the Jacobians of $\tilde{C},
C$ respectively and let $P = P(\tilde{C},\sigma) =
(1-\sigma )\tilde{J}$ be the Prym variety. One has the maps $\pi ^* :
J\longrightarrow \tilde{J} , Nm : \tilde{J}\longrightarrow J$ such
that $\pi^* \circ Nm = 1+\sigma $. We denote by $j : P
\longrightarrow  \tilde{J} $ the embedding. The kernel of $\pi ^* : J
\longrightarrow  \tilde{J} $ is $\{0,\eta \}$ where $\eta \in J_2$.

Conversely, given $C$ and $\eta \in J(C)_2 , \eta \neq 0$ one
constructs a double unramified covering $\pi :  \tilde{C}
\longrightarrow C$ such that $\pi _*{\cal O}_{\tilde{C}} \simeq {\cal
O}_C \oplus{\cal O}_C(\eta)$ and gets the set-up above.

    Let $J \simeq {\bf C}^g/\Lambda , \tilde{J} \simeq {\bf
C}^{2g-1}/\tilde{\Lambda} , P\simeq {\bf C}^{g-1}/\Lambda_{\_}$.
Choosing as in \cite{fay}, \cite{clem} symplectic bases
$\{a_0,...,a_{g-1},b_0 ,...,b_{g-1}\} ,
\{\tilde{a}_0,...,\tilde{a}_{2g-2},\tilde{b}_0,...,\tilde{b}_{2g-2}\}
, \{
\tilde{a}_1-\tilde{a}_g,...,\tilde{a}_{g-1}-\tilde{a}_{2g-2},\tilde{b}
_1-\tilde{b}_g,...,\tilde{b}_{g-1}-\tilde{b}_{2g-2}\} $
of $\Lambda , \tilde{\Lambda }$ and $\Lambda_{\_}$ respectively we
can assume that $\eta = \frac{1}{2}a_0(mod \Lambda)$. Let us denote
by $\tau \in {\cal H}_g , \tilde{\tau} \in {\cal H}_{2g-1}, \Pi \in
{\cal H}_{g-1}$ the corresponding period matrices. One has the
following formulas :
\begin{equation}\label{e8.1}
\begin{array}{ccc}
(\pi ^*)_*\left[ \begin{array}{cc}\alpha _0&\alpha \\\beta _0&\beta
\end{array}\right] & = & \left[ \begin{array}{rcc}\alpha _0&\alpha
&\alpha \\2\beta _0&\beta &\beta  \end{array}\right] \\\mbox{}\\
Nm_*\left[ \begin{array}{ccc}\alpha _0&\alpha &\alpha '\\\beta
_0&\beta &\beta ' \end{array}\right] & = & \left[ \begin{array}{rc}
2\alpha _0&\alpha +\alpha '\\\beta _0&\beta +\beta '
\end{array}\right] \\ \mbox{}\\ j_*\left[ \begin{array}{c}\alpha
\\\beta  \end{array}\right] &=&\left[ \begin{array}{ccc}0&\alpha
&-\alpha \\0&\beta &-\beta  \end{array} \right] \end{array}
\end{equation}
Here $\alpha ,\alpha ',\beta ,\beta ' \in  {\bf R}^{g-1} ,\alpha
_0,\beta _0\in {\bf R}$ and $(\pi ^*)_*, Nm_* , j_*$ are the
linear maps which induce the homomorphisms   $\pi^* , Nm , j$.

   Wirtinger's theorem \cite{fay} states that there is a symmetric
theta divisor $\Theta _0$ of $\tilde{J}$ which is equal to
$W_{\tilde{g}-1}-\pi ^*\Delta $ for a certain theta characteristic
$\Delta $ of $C$, such that $\Theta_0\mid _P = 2\Xi $ for some
symmetric
theta divisor $\Xi$ of  $P$. Moreover any point $L-\pi ^*\Delta $ of
$\Theta _0\cap P$ satisfies the properties
$$
Nm(L) \equiv K_C , h^0(\tilde{C},L) \equiv 0(mod 2) , h^0(\tilde{C},
L) \geq 2
$$
$\Theta _0$ is the  divisor of the theta function $\theta[\lambda ]
(z,\tilde{\tau})$, where
$$
\lambda = \frac{1}{2}\tilde{ a}_0 =
\left[ \begin{array}{ccc}0&0&0 \\\frac{1}{2}&0&0 \end{array}\right]
$$
Next we recall and state in more general form some results of Welters
\cite{wel}. Let $f : X \longrightarrow Y$ be a double covering of
smooth, projective curves. It might be ramified. Let $\Lambda = \mid
D \mid$ be a complete linear system on $Y$ and let $deg(D) = d$. One
denotes by $S$ the subscheme of $X^{(d)}$ which is the pull-back of
$\Lambda $ under $f^{(d)}$
$$
\begin{array}{ccl}
          S&\longrightarrow & X^{(d)}\\
\downarrow &                &\downarrow \\
   \Lambda &\longrightarrow & Y^{(d)}

\end{array}
$$
The arguments on pp.103-107 of \cite{wel} combined with Riemann-Roch's
 theorem give the following proposition
\begin{prop}\label{p91.1}
Let $\hat{D}$ be a closed point of $S$ and let $A$ be the maximal
effective divisor of $C$ such that $\hat{D}=\pi ^*A + E$ with $A\geq
0,E\geq 0$. Let $D=f^{(d)}(\hat{D})=2A+E_1$. Then

     (i) $S$ is nonsingular at $\hat{D}$ if and only if
$$
h^0(D-A) = h^0(D) - deg(A)
$$

   (ii) Suppose $S$ is nonsingular at $\hat{D}$. Then
$f^{(d)}\mid _{S} : S\longrightarrow \Lambda $ is nondegenerate at
$\hat{D}$ if and only if $f^{(d)}:X^{(d)}\longrightarrow Y^{(d)}$ is
nondegenerate at $\hat{D}$ and this is the case if and only if $D$
contains no branch points of $f$ and $A=0$.
\end{prop}
Let $B\subset Y$ be the branch locus of $f$ and let $\delta$ be the
invertible sheaf with $\delta ^{\otimes 2}\simeq {\cal O}_Y(B)$ which
determines the covering by $f_*{\cal O}_X\simeq {\cal O}_Y\oplus
\delta $. An effective divisor $E$ of $X$ is called $f$-simple if
$E\not \geq f^*(y)$ for any $y\in Y$. The following lemma is due to
Mumford \cite{mum}.
\begin{lem}\label{l92.1}
Let $A$ be a divisor of $Y$ and let $E$ be an effective $f$-simple
divisor of $X$. Then there is an exact sequence
$$
0\longrightarrow {\cal O}_Y(A)\longrightarrow f_*{\cal O}_X(\pi^*A+E)
\longrightarrow {\cal O}_Y(A+Nm_f(E))\otimes \delta
^{-1}\longrightarrow 0
$$
\end{lem}
\begin{cor}\label{c92.2}
Under the assumptions of the preceding lemma suppose that \\
\noindent $deg(A)+deg(E)<deg(\delta )$. Then
$$
h^0(X,\pi ^*A+E) = h^0(Y,A)
$$
\end{cor}

\section{Double unramified coverings of elliptic curves}\label{s2}
Let $\tilde{E}$ be an elliptic curve. Choosing a point $o\in
\tilde{E}$ we shall sometimes identify $\tilde{E}$ with $J(\tilde{E})$
by the map $x\mapsto cl(x-o)$. Let $\sigma :\tilde{E}\longrightarrow
\tilde{E}$ be an involution without fixed points.
\begin{lem}\label{l10.1}
There exists $\mu \in J(\tilde{E})$ of order $2$ such that $\sigma
(x) = x+\mu $. Furthermore $Ker(\pi_*:J(\tilde{E})\longrightarrow
J(E))=\{ 0,\mu \} $.
\end{lem}
{\bf Proof.} Let $\mu =\sigma (0)-0.$ Since $P(\tilde{E},\sigma
)=(\sigma -1)J(\tilde{E})=0$ we have $\sigma (x-o)\equiv x-o$, thus
$\sigma x=x+\mu.$ Furthermore $2(\sigma (o)-o)=(1-\sigma )(\sigma
(o)-0)\equiv 0$ and $\pi _*(\sigma (o)-o)\equiv 0$, thus
$Ker\pi _*=\{ 0,\mu \} $ since $\# Ker\pi _*=2.$ q.e.d.

Using the notation of the Introduction let $\tilde{E}\simeq {\bf
C}/\Lambda _{\tau }$ where $\Lambda _{\tau }={\bf Z}\tau +{\bf Z}$
with $Im(\tau )>0$ and $\mu = \frac{1}{2}(\tau +1)(mod \Lambda
_{\tau })$. Consider $\lambda (\tau )$ defined by Eq.~(\ref{ii.1}).
Then $\{ \lambda (\tau ),1/\lambda (\tau )\} $ is the set of the roots
 of the equation $x^2-k(\tau )x+1=0$  where
$$
k(\tau )=-(\theta _{10}(0,\tau )^8+\theta _{01}(0,\tau )^8)/
\theta _{10}(0,\tau )^4\theta _{01}(0,\tau )^4
$$
\begin{prop}\label{p12.1}
The map
\begin{equation}\label{e12.1}
k:\Gamma _{1,2}\backslash {\cal H}\longrightarrow {\bf C}-\{0,2\}
\end{equation}
is an analytic isomorphism.
\end{prop}
{\bf Proof.} Let $\Gamma _2$ be the level $2$  subgroup of $PSL(2,{\bf
 Z})$
$$
\Gamma _{2}=\{ \left( \begin{array}{cc}a&b\\c&d  \end{array}
 \right) \equiv \left( \begin{array}{cc}1&0\\0&1  \end{array}
\right) (mod 2)\}
$$
Then $\mid \Gamma _{1,2}:\Gamma _{2}\mid =2$ and the element
$S=\left( \begin{array}{cc}0&1\\-1&0  \end{array}\right) $,
 belongs to $\Gamma _{1,2}\backslash \Gamma _{2}$. It is well-known
(see e.g. \cite{clem}) that the map
$$
\lambda :\Gamma _{2}\backslash {\cal H}\longrightarrow {\bf C}-\{
0,1\}
$$
given by Eq.~(\ref{ii.1}) is an isomorphism. We have $S(\tau
)=-1/\tau $ and $\lambda (-1/\tau )=1/\lambda (\tau )$. The factor of
$\Gamma _{2}\backslash {\cal H}$ by the action of $S$ is $\Gamma
_{1,2}\backslash {\cal H}$, thus $k$ is an analytic isomorphism.
q.e.d.

Explicitly, given $k\neq 0,2$ we find $\lambda $ such that
$\lambda +1/\lambda =k$ and the corresponding pair $(\tilde{E},\mu \in
J(\tilde{E})_2)$ is given by the equation $y^2=x(x-1)(x-\lambda )$ and
 the point $\mu =cl(p_1-p_2)$ where $p_1=(0,0),p_2=(1,0)$.

\section{Extended Prym data}\label{s3}
Let $\tilde{C},C$ etc. be as in Section~(\ref{s1}). Let $\Theta _{0}$
be
the divisor of the theta function \[
\theta \left[ \begin{array}{ccc}0&0&0\\\frac{1}{2}&0&0
\end{array}\right]
(z,\tilde{\tau })\] Let us denote by $q_{0}$ the quadratic form
$q_{\Theta _0} : \tilde{J}_2\longrightarrow {\bf Z}_2$ defined by
Eq.~(\ref{e5.1}). By Eq.~(\ref{e6.1}) and (\ref{e6.2})
one has
\begin{equation}\label{e14.1}
q_{0}\left( \left[ \begin{array}{cc}\alpha \\\beta   \end{array}
\right] \right)=4\alpha ^t\beta +2\alpha _0(mod 2)
\end{equation}
Hence $q_{0}(\rho )=0$ for any $\rho \in P_2$. This follows also from
Wirtinger's theorem. The same property holds for any symmetric theta
divisor of the orbit $\{ \Theta _{0}+\rho \mid \rho \in P_2\} $. Let
us denote $\pi ^*(J)$ by $B$. By Eq.~(\ref{e8.1}) one has
$B_2\supset P_2$.

\begin{lem}\label{l14.1}
A symmetric theta divisor $\Theta \subset \tilde{J}$ has the property
that $q_{\Theta }$ vanishes on $P_2$ if and only if $\Theta =\Theta_
{0}+\alpha $ where $\alpha \in B_2$ and $q_{0}(\alpha )=0$.
\end{lem}
{\bf Proof.} Suppose $q_{\Theta }(P_2)=0$. Let $\Theta =\Theta_
{0}+\alpha $ for some $\alpha \in \tilde{J}_2$.Then for $\rho \in
P_2$ one has by Eq.~(\ref{e6.1})
$$
q_{\Theta }(\rho )=q_{0}(\rho +\alpha )=q_{0}(\rho )+e_2(\rho ,\alpha
)+q_{0}(\alpha )-q_{0}(0)
$$
Setting $\rho =0$ we conclude that $q_{0}(\alpha )=0$. Thus
$q_{\Theta }(P_2)=0$ implies that $e_2(P_2,\alpha )=0$.
Eq.~(\ref{e8.1}) show that the latter holds if and only if $\alpha
\in B_2$. Conversely, if $\alpha \in B_2$ and $q_{0}(\alpha )=0$,
then $q_{\Theta }(P_2)=0$ by the formula for $q_{\Theta }$ above.
q.e.d.

The zeros of $q_{0}$ which belong to $B_2$ are the following three
cosets with respect to $P_2$~: $P_2,\lambda _1+P_2,\lambda _2+P_2$
where $\lambda _1=\frac{1}{2}\tilde{a}_0(mod \tilde{\Lambda }),
\lambda _2=\frac{1}{2}\tilde{a}_0+\frac{1}{2}\tilde{b}_0(mod
\tilde{\Lambda })$. Let $\Theta _{1}=\Theta _{0}+\lambda _1,\Theta
_{2}=\Theta _{0}+\lambda _2$. We conclude that there are three $P_2$ -
orbits of
symmetric theta divisors $\Theta $ such that $q_{\Theta }$ vanishes on
 $P_2$ :
\[ \begin{array}{ccc}
\{ \Theta _{0}+\rho \}\; ,&\{ \Theta _{1}+\rho \}\; ,&\{ \Theta
_{2}+\rho \} \end{array}\]
These are respectively the divisors of the theta functions
\[ \begin{array}{ccc}
\theta \left[ \begin{array}{ccc}0&\alpha
&\alpha \\\frac{1}{2}&\beta &\beta \end{array}\right](z,\tilde{\tau
})\; ,& \theta \left[ \begin{array}{ccc}0&\alpha
&\alpha \\0&\beta &\beta \end{array}\right](z,\tilde{\tau })\; ,&
\theta \left[
\begin{array}{ccc}\frac{1}{2}&\alpha &\alpha \\0&\beta
&\beta \end{array}\right](z,\tilde{\tau })
\end{array}\]
where $\alpha ,\beta \in \frac{1}{2}{\bf Z}^{g-1}$.
\begin{lem}\label{16.2}
Let $\mu = \frac{1}{2}\tilde{b}_0(mod \tilde{\Lambda })$. Then
$Ker(Nm)=P\cup P_{\mu } , Ker(Nm)\cap B_2=P_2\cup (\mu +P_2)$ and for
any $x\in \tilde{C}$ there exists a unique $\xi \in P$ such that
$\sigma x-x\equiv \mu +\xi $.
\end{lem}
{\bf Proof.} It is well-known \cite{fay},\cite{mum} that $Ker(Nm)$
 has two connected components $P$ and $P_{\_ }$. Since $\mu \in
B_2\backslash P_2$ and $Nm(\mu )=0$ we get that $P_{\_ }=P_{\mu }$.
Hence $(P\cup P_{\_ })\cap B_2=P_2\cup (P_2+\mu )$. The last
statement of the lemma follows from the equality $\sigma x-x+P=
P_{\_ }$ \cite{mum}.
 q.e.d.

Let $\tilde{J}_{2g-2},J_{2g-2}$ be the divisor classes of degree
$2g-2$ on $\tilde{C},C$ respectively and let
$Nm:\tilde{J}_{2g-2}\longrightarrow J_{2g-2}$ be the norm map. The
subvariety $Nm^{-1}(K_C+\eta )$ is a principal homogeneous space for
$Nm^{-1}(0)=P\cup P_{\_ }$, thus it has two connected components. Let
\begin{equation}\label{e165.1}
Z=\{ L\in \tilde{J}_{2g-2}\mid Nm(L)=K_C+\eta ,\; h^0(L)\geq 1\}
\end{equation}
and let $Z=Z_1\cup Z_2$ where $Z_{i}$ are the intersections of $Z$
with the connected components of $Nm^{-1}(K_C+\eta )$.

\begin{lem}\label{l1634.1}
Let $M$ be an effective divisor of $\tilde{C}$ such that $Nm(M) \in
\mid K_C+\eta \mid$. Suppose $x\in \tilde{C}$ and $x\not \in Bs\mid
M\mid $. Then
$$
h^0(M+\sigma x-x)=h^0(M)-1
$$
\end{lem}
{\bf Proof.} By Riemann-Roch's theorem $x$ is a base point of
$\mid K_{\tilde{C}}-M+x\mid $. Now, $K_{\tilde{C}}-M\equiv\sigma (M)$
, thus $\sigma (x)$ is a base point of $\mid M+\sigma (x)\mid $ which
proves the lemma since $x\not \in Bs\mid M\mid $. q.e.d.

\begin{prop}\label{p17.1}
The theta divisors $\Theta _{1}$ and $\Theta _{2}$ do not contain $P$
. The restrictions $\Theta _i.P$ are connected and reduced divisors of
 $P$ which belong to the linear system $\mid 2\Xi \mid $. Furthermore,
 up to a possible reordering of $Z_i$ one has
$$
\Theta _i.P=Z_i-\pi ^*\Delta -\lambda _i
$$
The point $L-(\pi ^*\Delta +\lambda _i)$ is nonsingular if and only if
 $h^0(\tilde{C},L)=1$. The corresponding tangent hyperplane is equal
to $Nm(\mid L\mid )\in \mid K_C\otimes \eta \mid $ via the
identification
$$
T_0(P)^*\simeq H^0(\tilde{C},K_{\tilde{C}})^-\simeq H^0(C,K_C\otimes
\eta )
$$
\end{prop}
{\bf Proof.} Let $\kappa _i=\pi ^*\Delta +\lambda _i,i=1,2$.We have
$Nm(\lambda _i)=\eta $, thus $Nm(\kappa _i)=K_C+\eta $ and for
$L$ with $h^0(L)\geq 1$ one has $L-\kappa _i\in\Theta _i\cap Ker(Nm)$
if and only if $Nm(\mid L\mid )\in \mid K_C+\eta \mid $. Since
$dim\mid K_C+\eta \mid =g-2$ we conclude that neither $P$ nor $P_{\_
}$ are contained in $\Theta _i$. Furthermore $\Theta _i$ are ample,
so $\Theta _i\cap P$ are not empty and $\Theta _i.P$ are connected
divisors of $P$. Upon a possible reordering of $Z_1$ and $Z_2$ we
have $Z_1-\kappa _1=\Theta _1\cap P,Z_2-\kappa _1=\Theta
_1\cap P_{\_ }$. Since $P_{\_ }=P+\mu $ and $\lambda _2=\lambda
_1+\mu $ we obtain $Z_2-\kappa _2=\Theta _2\cap P$.

   {\bf Claim.} {\it For every irreducible component $T$  of any
$\Theta _i.P$ the general element $L-\kappa _i\in T$ satisfies
$h^0(\tilde{C},L)=1$ }.\\
{\bf Proof.} Suppose the contrary. Then for any $x\in \tilde{C}$ the
image of the map
\[ \begin {array}{cc}
\psi :T\times \tilde{C}\longrightarrow Ker(Nm)\; ,&
\psi (L,x)=L+\sigma x-x
\end{array}\]
is contained in $Z$. This image must be of dimension $dimT+1=g-1$.
Indeed, if $M=L+\sigma x-x$, then for every sufficiently general
$L\in T$ and $x\in \tilde{C}$ one has by Lemma~(\ref{l1634.1}) that
$h^0(M)=h^0(L)-1$. If $dim\psi (T\times \tilde{C})\leq dimT$, then
for every sufficiently general $M\in Im(\psi ),x\in \tilde{C}$ one
has $h^0(M-\sigma x+x)=h^0(L)=h^0(M)+1$ which is an absurd by
Lemma~(\ref{l1634.1}). Now, $dimZ=g-2$, thus it is impossible that
$dim\psi (T\times \tilde{C})=g-1$. q.e.d.

Now, suppose that $L-\kappa _i\in\Theta _i\cap P$ is an element with
$h^0(L)=1$ and let $D=\mid L\mid $. Since $Nm(D)\in\mid K_C+\eta\mid
$ there is an anti invariant holomorphic differential $\omega $ of
$\tilde{C}$ whose divisor of zeros is $\pi ^*D$. Since $h^0(L)=1$
the point $L-\kappa _i$ is nonsingular on $\Theta _i$ and the tangent
space in $T_0\tilde{J}$ is given by the equation $\omega =0$. We see
that $\Theta _i$ and $P$ intersect transversely at $L-\kappa _i$ and
 the tangent hyperplane of $\Theta _i.P$ at $L-\kappa _i$ is given by
the same equation $\omega =0$ in $T_0P$ since $\omega $
is anti invariant. What we have proved implies also that $Sing(\Theta
_i.P)=P\cap Sing\Theta _i$. This concludes the proof of the
proposition. q.e.d.

\begin{cor}\label{c19.1}
All irreducible components of $Z$ are of dimension $g-2$.
\end{cor}
We see that the $P_2$-orbit $\{ \Theta _0+\rho  \} $ is distinguished
among the three $P_2$-orbits of symmetric theta divisors $\Theta$
which satisfy $q_{\Theta}\mid _{P_2}=0$ by the property that the
restriction of any $\Theta _0+\rho $ on $P$ is equal to
twice a theta divisor of $P$. It is also distinguished by the
property that every $\{ \Theta _0+\rho  \} $ contains $P_{\_ }=P+\mu $
. Indeed, the fact that $\Theta _0\supset P+\mu $ follows from the
parity lemma \cite{tju} and is well-known \cite{mum1},\cite{fay}. If
$\Theta _1+\rho $ contained $P_{\_ }=P+\mu $, then $\Theta _2$
would contain $P$ since $\Theta _2=\Theta _1+\mu $ and $\rho  \in P_2$
which contradicts Proposition~(\ref{p17.1}). Notice that the latter
distinction of $\{ \Theta _1+\rho  \} $ and $\{ \Theta _2+\rho  \} $
parallels the distinction of $\theta _{10}(z,\tau )$ and $\theta
_{01}(z,\tau )$ in the elliptic case. We can now state our extension
of the Prym data.

{\sc Extended Prym data}.{\it One associates to every
algebraic, smooth, irreducible, projective curve $\tilde{C}$ of genus
$\geq 3$ with an involution $\sigma
:\tilde{C}\longrightarrow \tilde{C}$ without fixed points, the
principally polarized Prym variety $(P,\Xi )$ and the two $P_2$-orbits
 $O_1,O_2\subset \mid 2\Xi \mid $ which consist of the restrictions
$\Theta .P$ of the symmetric theta divisors $\Theta \subset
J(\tilde{C})$ such that $q_{\Theta }\mid _{P_2}=0$ and $\Theta \not
\supset P_{\_ }$ where $Ker(Nm_{\pi }:J(\tilde{C})\longrightarrow
J(C))=P\cup P_{\_ }$ }

We can now state our result which is a kind of generalization of
Proposition~(\ref{p12.1}) to curves of genus $>1$.

\begin{theo}\label{t20.1}
The pair $(\tilde{C},\sigma )$ is uniquely determined up to
isomorphism
by the extended Prym data $(P(\tilde{C},\sigma ),\Xi ),O_1,O_2\subset
\mid 2\Xi \mid $.
\end{theo}

\section{The semicanonical map and the Gauss map}\label{s4}
In this section $K_C\in Pic^{2g-2}(C)$ is the canonical sheaf of $C$
and
$\eta \in Pic^0(C)_2$ is the sheaf with $\eta ^{\otimes 2}\simeq
{\cal O}_C$ such that $\pi _*{\cal O}_{\tilde{C}}\simeq {\cal
O}_{C}\oplus \eta $. We shall denote by $\varphi  _K,\varphi
_{K\otimes
\eta }$ the canonical, respectively semicanonical map of the curves
under consideration. Let $L=K_C\otimes \eta $. The following lemma
follows elementary from Riemann-Roch's theorem.

\begin{lem}\label{l21.1}
Suppose $g(C)\geq 2$. Then $\mid K_C\otimes \eta \mid $ has base
points if and only if $C$ is hyperelliptic and $\eta \simeq {\cal
O}_{C}(p_1-p_2)$ where $p_1,p_2$ are ramification points for the
double covering $f:C\longrightarrow {\bf P}^1$. In this case
$p_1+p_2=Bs\mid K_C\otimes \eta \mid $ and $K_C\otimes \eta \simeq
(f^*{\cal O}_{{\bf P}^1}(g-2))(p_1+p_2)$.
\end{lem}

\begin{lem}\label{l21.2}
Suppose $g(C)\geq 3 $ and $\mid K_C\otimes \eta \mid $  is without
base points. Then \\
$\varphi  _L:C\longrightarrow \mid K_C\otimes \eta \mid
^*={\bf P}^{g-2}$ is a birational embedding except in the following
two cases :

(i) $g(C)=3$. Then $\varphi  _L:C\longrightarrow {\bf P}^1$ is
of
degree $4$.

(ii) $g(C)\geq 4$, $C$ is bi-elliptic, i.e. it is a double
covering $f:C\longrightarrow E$ of an elliptic curve, and $\eta \simeq
f^*(\epsilon )$ where $\epsilon \in Pic^0(E)_2$. Here $\varphi
_L=\varphi
_{\delta \otimes \epsilon }\circ f$ where $\delta $ is the invertible
sheaf of $E$  which determines the covering, i.e. $\delta ^{\otimes
2}\simeq {\cal O}_{E}(x_1+...+x_{2g-2})$ for the branch points
$x_1,...,x_{2g-2}$ and $f_*{\cal O}_{C}\simeq {\cal O}_{E}\oplus
\delta $.
\end{lem}
{\bf Proof.} The case $g(C) =3$ is clear, so let us suppose that
$g\geq 4$. Let $X=\varphi  _L(C)$ and let $d$ be the degree of the
map
$\varphi  _L:C\longrightarrow X$. We have $d.deg(X)=2g-2$ and
$deg(X)\geq g-2$. This implies that the case $d\geq 3$ may occur only
if $g=4,deg(X)=2,d=3$. Otherwise either $\varphi  _L$ is a birational
embedding or $d=2,deg(X)=g-1 $. In the latter case $p_a(X)=1$.
Suppose $d=2$ and $X$ is singular. Then the normalization of $X$ is
$\hat{X}\simeq {\bf P}^1$ and we can decompose $\varphi  _L$ as
$$
\varphi  _L=g\circ f:C\longrightarrow \hat{X}\longrightarrow {\bf
P}^{g-2}
$$
Since $\varphi  _L$ is obtained from a complete linear system, $g$
must
have the same property, thus $g^*{\cal O}_{{\bf P}^{g-2}}(1)\simeq
{\cal O}_{{\bf P}^1}(g-2)$. This is impossible since $deg(X)=g-1$.
Consequently if $d=2$ then $X\subset {\bf P}^{g-2}$ is an elliptic
curve. Let $E=X,f=\varphi  _L:C\longrightarrow E.$  Since $K_C\simeq
f^*(\delta )$ and $K_C\otimes \eta \simeq f^*{\cal O}_{E}(1)$ we
conclude that $\eta \simeq f^*(\epsilon )$ for some $\epsilon
\in Pic^0(E).$ The covering $f$
is ramified, so $f^*:Pic^0(E)\longrightarrow Pic^0(C)$ is an
injection, hence $\epsilon \in Pic^0(E)_2$
and $\epsilon \not \simeq {\cal O}_{E}$. Conversely, if
$f:C\longrightarrow E$ is a double covering of an elliptic curve, and
$\eta =f^*(\epsilon )$ with $\epsilon \in Pic^0(E)_2,\epsilon \not
\simeq {\cal O}_{E}$ then
$$
H^0(C,K_C\otimes \eta)\simeq H^0(C,\pi ^*(\delta \otimes \epsilon
))=\pi ^*H^0(E,\delta \otimes \epsilon )
$$
Hence $\varphi  _L=\varphi  _{\delta \otimes \epsilon }\circ f$ and
$d=2$.

It remains to rule out the possibility $g=4,d=3,deg(X)=2$. Here \\
$f=\varphi  _L:C\longrightarrow X\simeq {\bf P}^1$ so $L\simeq
M^{\otimes
2}$, where $deg(M)=3,h^0(C,M)=2$, and $\mid M\mid $ is without base
points. Thus $C$ is not hyperelliptic and
$$
\varphi  _K(C)=Q\cap F\subset {\bf P}^3
$$
where $Q$ is a quadric and $F$ is a cubic surface. We have $\mid L\mid
=\mid M\mid +\mid M\mid $ since $dim\mid L\mid =2$. Let $l_1,l_2$ be
lines in $Q$ such that $l_1+l_2=Q.H$ for a plane $H$. Let $M_i={\cal
O}_{C}(C.l_i).$ We can assume that $M=M_1$. Then $\eta \simeq
M_1\otimes M_2^{-1}.$ If $Q$ were singular, then $M_1\simeq M_2,$ so
$\eta \simeq {\cal O}_{C}$ which is absurd. Suppose $Q$ is
nonsingular.
Then $\eta ^{\otimes 2}\simeq {\cal O}_{C}$ implies $\mid
M_2^{\otimes 2}\mid =\mid M_1^{\otimes 2}\mid =\mid M_1\mid +\mid
M_1\mid $. This is again impossible since any reduced divisor
$x_1+x_2+x_3\in \mid M_2\mid $ can have only two common points with
any two divisors $D_1,D_2\in \mid M_1\mid $. q.e.d.

Suppose $g\geq 3$. Following Welters \cite{wel} let $S$ be the
subscheme of $\tilde{C}^{2g-2}$ which is the pull-back of $\mid
K_C\otimes \eta \mid \subset C^{2g-2}$
$$
\begin{array}{ccl}S&\longrightarrow &\tilde{C}^{2g-2}\\
\downarrow &\mbox{}&\downarrow Nm\\\mid K_C\otimes \eta \mid
&\longrightarrow &C^{2g-2} \end{array}
$$
It breaks naturally into two disjoint subschemes $S=S_1\cup S_2.$ The
singularities of $S$ can be calculated by Proposition~(\ref{p91.1})
 with
 $X=\tilde{C},Y=C,f=\pi ,\pi ^{(2g-2)}=Nm$. Since $S$ is a locally
complete intersection and $\pi ^{(2g-2)}$ is a finite map every
irreducible component of $S$ has dimension $g-2$. A Zariski open,
dense subset of $\mid K_C\otimes \eta \mid $ consists of reduced
divisors by Lemmas~(\ref{l21.1}) and (\ref{l21.2}), so according
to
Proposition~(\ref{p91.1})  $S$ is reduced. The subvarieties
$S_1,S_2$ are connected provided $g\geq 3$ \cite{wel}.

Let $T_1,T_2$ be divisors from $\mid 2\Xi \mid $ which belong to the
orbits $O_1,O_2$ respectively. Suppose $g\geq 3$.The Gauss maps
$G_i:T_i^{ns}\longrightarrow {\bf P}(T_0P)^*$ are defined on the
nonsingular loci of $T_i$ and send a point $x\in T_i^{ns}$ to the
translation of the tangent hyperplane $T_x(T_i)$ to $0\in P$. Let
$T=T_1\sqcup T_2$ and let $G:T^{ns}\longrightarrow {\bf P}(T_0P)^*$
 be the
 map whose restriction on $T_i^{ns}$ equals $G_i$. Let $S^0$ be the
Zariski open subset of $S$ which consists of those $\hat{D}$ with
$h^0(\tilde{C},\hat{D})=1$. By Proposition~(\ref{p17.1}) the map
$cl:S^0\longrightarrow Z^{ns}$ is an isomorphism and moreover one can
identify $T_i^{ns}$ with $Z_i^{ns}$ by translation and ${\bf
P}(T_0P)^*$ with $\mid K_C\otimes \eta \mid $. Since the Gauss map
does not depend on the translation one has for the Gauss map
$G:Z^{ns}\longrightarrow \mid K_C\otimes \eta \mid $ and every
$\hat{D}\in S^0$ the formula
\begin{equation}\label{e26.1}
G(cl(\hat{D}))=Nm(\hat{D})
\end{equation}
\begin{prop}\label{p26.1}
Let $L=K_C\otimes \eta $. Suppose $g\geq 4$. Let $R\subset T_i$ be the
 ramification locus of $G$ and let $B$ be the algebraic closure of
$G(R)$.

(i) If $C$ is hyperelliptic and $\eta \simeq {\cal
O}_C(p_1-p_2)$, where $p_1,p_2$ are Weierstrass points of $C$,
then $P(\tilde{C},\sigma )\simeq J(C_2)$ for a certain hyperelliptic
curve $C_2$ (see Section~(\ref{s5})) and
$$
B=\varphi  _K(C_2)^*\cup \bigcup _{i=1}^{2g}\varphi  _K(q_i)^*
$$
where $\varphi  _K(C_2)^*$ is the dual hypersurface of $\varphi
_K(C_2)$
and $\varphi  _K(q_i)^*$ are the stars of hyperplanes which contain
$\varphi
_K(q_i)$, where $q_i$ are the Weierstrass points of $C_2$.

(ii) If $\mid K_C\otimes \eta \mid $ is without base points and
$\varphi  _L:C\longrightarrow \mid K_C\otimes \eta \mid ^*$ is a
birational embedding, then $B$ has a unique irreducible component of
dimension $g-3$ and degree $>1$. This component is equal to
$\varphi  _L(C)^*$.

(iii) If $\mid K_C\otimes \eta \mid $ is without base points and
$$
f=\varphi  _L:C\longrightarrow  \varphi  _L(C)=E
$$
is a map of degree $2$ onto an elliptic curve, then
$$
B=E^*\cup \bigcup_{i=1}^{2g-2}x_i^*
$$
where $x_i$ are the branch points of $\varphi  _L$.
\end{prop}
{\bf Proof.} (i) In this case any $T_i$ is the union of two translates
 of the theta divisor $\Xi \subset P$ as it will be shown in
Section~(\ref{s5}). So, Part(i) follows from the description of the
branch locus of the Gauss map of the theta divisor of a hyperelliptic
Jacobian \cite{and}.

Now, let us assume that $\mid K_C\otimes \eta \mid $ is without base
points. If $Nm:S^0\longrightarrow \mid K_C\otimes \eta \mid $ is
degenerate at $\hat{D}\in S^0$, then by Proposition~(\ref{p91.1})
$\hat{D}=\pi ^*A+E$ for some $A>0$, so $D=Nm(\hat{D})=2~A+Nm(E)$. Let
$H\subset \mid K_C\otimes \eta \mid ^*$ be the hyperplane which
corresponds to $D$. Either $H$ contains an image of a ramification
point of the map $\varphi  _L:C\longrightarrow \mid K_C\otimes \eta
\mid
^*$, or $H$ is tangent to a branch $\varphi  _L(U)$ at a point
$\varphi
_L(p)$, where $p\in U\subset C$ and $\varphi  _L$ is nondegenerate at
$p$.
 The former case can happen only for finitely many points. This shows
that any component of $B$ of dimension $g-3$ must be either a star of
hyperplanes which contain a branch point $\varphi  _L(p)$, or it is
contained in the dual variety $\varphi  _L(C)^*$.

{\bf Proof of (ii).} It remains to show that $\varphi  _L(C)^*\subset
B$.
Let the hyperplane $H\subset \mid K_C\otimes \eta \mid ^*$ be a
sufficiently general element of $\varphi  _L(C)^*$ and let
$$
D=\varphi  _L^*(H)=2p+p_3+...+p_{2g-2}
$$
be the corresponding divisor of $\mid K_C\otimes \eta \mid $. Here
$p,p_3,...,p_{2g-2}$ are distinct points of $C$ and $\varphi  _L$ is
not
degenerate at $p$.

{\bf Claim 1.} {\it  Let $\pi ^{-1}(p)=p'+p''$ One can choose
$p'_i\in
\tilde{C}$ with $\pi (p'_i)=p_i$ so that
$$
\hat{D}=p'+p''+p'_3+...+p''_{2g-2}
$$
has the property $h^0(\hat{D})=1$.}\\
{\bf Proof.} Let us choose arbitrary points $q_i\in \tilde{C}$ such
that $\pi (q_i)=p_i$. Let
$$
\hat{D}_0=p'+p''+q_3+...+q_{2g-2}
$$
If $h^0(\hat{D})\geq 2$, then at least one of the points $q_i$ is not
a base point of $\mid \hat{D}_0\mid $. Indeed, otherwise
$$
h^0({\cal O}_{\tilde{C}}(\pi ^*p))=h^0({\cal O}_{C}(p))+h^0({\cal
O}_{C}(p)\otimes \eta )\geq 2
$$
which is possible only in the Case (i). If $q_i\not \in Bs\mid
\hat{D}_0\mid $, let $\hat{D}_1=\hat{D}_0+\sigma (q_i)-q_i$. Then by
Lemma~(\ref{l1634.1}) $h^0(\hat{D}_1)=h^0(\hat{D}_0)-1$. Repeating the
 same argument with $\hat{D}_1$ etc. we obtain eventually the required
 divisor $\hat{D}$. q.e.d.

   Now, $\hat{D}\in S^0$, the Gauss map at the point $cl(\hat{D})\in
Z^{ns}$ equals $H$ by Eq.~(\ref{e26.1}) and it is ramified at
$\hat{D}$ according to Proposition~(\ref{p91.1}). So, $\varphi
_L(C)^*\subset B$.\\
   {\bf Proof of (iii).} We have proved above that
$$
B\subset E^*\cup \bigcup_{i=1}^{2g-2}x_i^*
$$
Let $H$ be a sufficiently general element of $E^*$ and let
$$
D=\varphi  _L^*(H)=2p+2q+p_5+...+p_{2g-2}
$$
be the corresponding divisor of $\mid K_C\otimes \eta \mid $. Here
$p,q,p_5,...p_{2g-2}$ are distinct points of $C,\varphi  _L$ is
nondegenerate at $p,q$ and $\{ p,q\} =f^*(x)$
for some $x\in E$.

    {\bf Claim 2.} {\it Let $\pi ^{-1}(p)=\{ p',p''\} ,\pi ^{-1}(q)=\{
q',q''\} $. One can choose $p'_i\in \tilde{C}$ with $\pi (p'_i)=p_i$
so that
$$
\hat{D}=p'+p''+2q'+p'_5+...+p'_{2g-2}
$$
has the property $h^0(\hat{D})=1$.}\\
   {\bf Proof.} Let $\eta = f^*(\epsilon )$ and let $y\in E $
be the point such that $\epsilon \simeq {\cal O}_{E}(y-x)$. Then one
has canonical isomorphisms
$$
\begin{array}{l}
H^0({\cal O}_{\tilde{C}}(\pi ^*(p+q))\simeq \pi ^*(H^0({\cal
O}_{C}(p+q))\oplus H^0({\cal O}_{C}(p+q)\otimes \eta))\\
\simeq (\pi ^*\circ f^*)(H^0({\cal O}_{E}(x))\oplus H^0({\cal
O}_{E}(y)))\simeq {\bf C}^2
\end{array}
$$
Thus $\mid p'+p''+q'+q''\mid $ is a pencil without base points. Using
the same argument as in Claim~1 we conclude that one can choose $p'_i$
 so that
$$
\hat{D}'=p'+p''+q'+q''+p'_5+...+p'_{2g-2}
$$
has the property $h^0(\hat{D}')=2$. Applying once more
Lemma~(\ref{l1634.1}) we conclude that $\hat{D}=\hat{D}'+ q'
-q''$ satisfies $h^0(\hat{D})=1$. q.e.d.

    We conclude as in Part (ii) that $E^*\subset B$. If $x_i=\varphi
_L(q_i)$ is a branch point of $f$ and $H$ is a sufficiently general
hyperplane in $\mid K_C\otimes \eta \mid ^*$ which contains $x_i$,
 then the corresponding divisor of $\mid K_C\otimes \eta \mid $ has
the form
$$
D=\varphi  _L^*(H)=2q_i+p_3+...+p_{2g-2}
$$
where $q_i,p_3,...,p_{2g-2}$ are distinct. The same argument as in
Part (ii) proves that $H\in B$, so $x_i^*\subset B$.
Proposition~(\ref{p26.1}) is proved. q.e.d.

We see that if $g\geq 4$ then the branch locus $B$ of the Gauss map
$G:T^{ns}\longrightarrow {\bf P}(T_0P)^*$ has a unique irreducible
component $B_0$ of dimension $g-3$ and degree $\geq 2$. We have shown
above that $B=X^*$ for a certain irreducible curve $X\subset {\bf
P}(T_0P)^*$. So, for $g\geq 4$, by the equality $(X^*)^*=X$
\cite{kleim} we obtain that $B_0^*$ is a curve
$X$. The following alternative takes place.
\begin{list}{(\roman{bean})}{\usecounter{bean}}
\item deg(X)=g-2
\item deg(X)=2g-2
\item deg(X)=g-1
\end{list}
The three cases correspond to those in Proposition~(\ref{p26.1}). If
Case (ii) occurs we prove Theorem~(\ref{t20.1}) as follows: $C$ is
isomorphic to the normalization of $X$. The normalization map
$f:C\longrightarrow X\subset {\bf P}(T_0P)^*$ is associated to the
complete linear system $\mid K_C\otimes \eta \mid $. Thus we obtain
$\eta \in J(C)_2$.

Cases (i) and (iii) are considered respectively in
Sections~(\ref{s5})
and (\ref{s6}). The case $g=3$ is treated in Section~(\ref{s7})
and the
case $g=2$ in Section~(\ref{s5}).

\section{The hyperelliptic case, $g\geq 2$}\label{s5}
Throughout this section we suppose that $C$ is a hyperelliptic curve
of genus $g\geq 2$,\\
$f:C\longrightarrow {\bf P}^1$ is the double covering
and $\eta \simeq {\cal O}_{C}(p_1-p_2)$ where $p_1,p_2$ are
ramification points of $f$ and $p_1\neq p_2$. Let $R=\{
p_1,p_2,p_3,...,p_{2g+2}\} $ be the set of ramification points of
$f,R_1=\{ p_1,p_2\}, R_2=R\backslash R_1$. Let $B_i=f(R_i),i=1,2$.
According to \cite{mum},\cite{dal} the covering \\$f\circ
\pi:\tilde{C}\longrightarrow C\longrightarrow {\bf P}^1$ has Galois
group ${\bf Z}_2\times {\bf Z}_2$. In the corresponding diagram of
Fig.1
\begin{figure}\label{f35.1}
\begin{center}
\begin{picture}(46,54)(0,0)
\put (19,0){\makebox(8,8){${\bf P}^1$}}
\put (0,23){\makebox(8,8){$C_1$}}
\put (19,23){\makebox(8,8){$C$}}
\put (38,23){\makebox(8,8){$C_2$}}
\put (19,46){\makebox(8,8){$\tilde{C}$}}
\put (19,46){\vector(-1,-1){15}}
\put (27,46){\vector(1,-1){15}}
\put (23,46){\vector(0,-1){15}}
\put (23,23){\vector(0,-1){15}}
\put (4,23){\vector(1,-1){15}}
\put (42,23){\vector(-1,-1){15}}
\put (24,15){$f$}
\put (24,38){$\pi $}
\put (6,38){$\pi _1$}
\put (37,38){$\pi _2$}
\put (6,15) {$f_1$}
\put (37,15) {$f_2$}
\end{picture}
\end{center}
\caption{}
\end{figure}
 $f_i:C_i\longrightarrow {\bf P}^1$ is branched at $B_i,i=1,2$.
Furthermore $\pi_2^*:J(C_2)\longrightarrow P(\tilde{C},\sigma )$ is an
 isomorphism. Let $\Theta _0=W_{\tilde{g}-1}(\tilde{C})-\pi ^*\Delta $
 be as in Section~(\ref{s3}).

\begin{lem}\label{l36.1}
Let $\Xi \subset P(\tilde{C},\sigma )=P$ be a symmetric theta divisor.
Then there exists a unique $\rho \in P_2$ such that

\begin{equation}\label{e36.1}
\Xi = \pi _1^*(\zeta _1)+\pi _2^*W_{g-2}(C_2)-\pi ^*\Delta-\rho
\end{equation}
where $\zeta _1$ is the rational equivalence class of the points of
the rational curve $C_1$.
\end{lem}
{\bf Proof.} By Wirtinger's theorem there is a unique translation
$\Theta =\Theta _0+\rho  $ with $\rho  \in P_2$ such that
$\Theta .P=2\Xi $. The points of $\Theta \cap P$ have the form
$L-\pi ^*\Delta-\rho  $
where $Nm(L)=K_C\; ,\; h^0(\tilde{C},L)\equiv 0(mod 2)$ and
$h^0(\tilde{C},L)\geq 2$. Now, $\mid K_C\mid \simeq f^*\mid {\cal
O}_{{\bf P}^1}(g-1)\mid $. One easily checks that if $\hat{D}$ is
effective divisor of $\tilde{C}$ and $Nm(\hat{D})\in f^*\mid
{\cal O}_{{\bf P}^1}(g-1)\mid $ then
$$
\hat{D}=\pi _1^*E+\pi _2^*F
$$
where $E$ and $F$ are effective divisors of $C_1,C_2$ respectively. We
 have $deg(E)+deg(F)=g-1$ which gives only two irreducible components
 of dimension $\geq g-2$ of
$$
Nm^{-1}(K_C)\cap W_{2g-2}(\tilde{C})
$$
namely $\pi _2^*W_{g-1}(C_2)$ of dimension $g-1$ and $\pi _1^*(\zeta
_1)+\pi _2^*W_{g-2}(C_2)$ of dimension $g-2$. On the other hand the
above intersection has, by the general theory, two irreducible
components: a translation of $\Xi $ and a translation of $P_{\_ }$.
This shows Eq.~(\ref{e36.1}). q.e.d.

Now, let us calculate the orbits $O_1,O_2\subset \mid 2\Xi \mid $. Let
 $T_i\in O_i,i=1,2$. By Proposition~(\ref{p17.1}) one has
\begin{equation}\label{e38.1}
T_i=Z_i-\pi ^*\Delta -\nu _i
\end{equation}
for some $\nu _i\in \lambda _i+P_2 , i=1,2.$
\begin{lem}\label{l38.1}
One can enumerate $\pi ^{-1}(p_1)$ as $\{ p'_1,p''_1\} $ and $\pi
^{-1}
(p_2)$ as $\{ p'_2,p''_2\} $ so that
\[ \begin{array}{ccc}Z_1&=&\pi _2^*W_{g-2}(C_2)+p'_1+p'_2\cup
\pi _2^*W_{g-2}(C_2)+p''_1+p''_2\\Z_2&=&\pi
_2^*W_{g-2}(C_2)+p'_1+p''_2\cup \pi _2^*W_{g-2}(C_2)+p''_1+p'_2
\end{array}\]
\end{lem}
{\bf Proof.} From Lemma~(\ref{l21.1}) we have
$$
\mid K_C\otimes \eta \mid =f^*\mid {\cal O}_{{\bf P}^1}(g-2)\mid
+p_1+p_2
$$
One easily checks that if $\hat{D}$ is effective divisor of
$\tilde{C}$ and $Nm(\hat{D})\in f^*\mid {\cal O}_{{\bf P}^1}(g-2)\mid
$ then
$$
\hat{D}=\pi _1^*E+\pi _2^*F
$$
where $E,F$ are effective divisors of $C_1,C_2$ respectively. We have
$deg(E)+deg(F)=g-2.$
Thus the only irreducible component of dimension $\geq g-2$ of
$$
Nm^{-1}(f^*{\cal O}_{{\bf P}^1}(g-2))\cap W_{2g-4}(\tilde{C})
$$
is $\pi _2^*W_{g-2}(C_2)$. The irreducible components of $Z=Z_1\cup
Z_2$ are of dimension $g-2$ by Corollary~(\ref{c19.1}) and the
transformation $L\mapsto L+\sigma (p)-p$ interchanges the two
components of $Nm^{-1}(K_C\otimes \eta)$. This shows that $Z_1$ and
$Z_2$ have the form given in the lemma. q.e.d.

Lemmas~(\ref{l36.1}),(\ref{l38.1}) and Eq.~(\ref{e38.1}) give the
following corollary
\begin{cor}\label{c39.1}
Let $\Xi $ be an arbitrary symmetric theta divisor of
$P(\tilde{C},\sigma )$ and let $T_1,T_2$ be two divisors of the orbits
 $O_1,O_2\subset \mid 2\Xi \mid $ respectively. Then
\[ \begin{array}{c}T_1=\Xi +p'_1+p'_2-\pi _1^*(\zeta _1)-\mu _1\cup
\Xi +p''_1+p''_2-\pi _1^*(\zeta _1)-\mu _1\\T_2=\Xi +p'_1+p''_2-\pi
_1^*(\zeta _1)-\mu _2\cup\Xi +p''_1+p'_2-\pi _1^*(\zeta _1)-\mu _2
\end{array}\]
for some $\mu _i\in \lambda _i+P_2,i=1,2$.
\end{cor}

Let us choose in an arbitrary way $\Xi ,T_1,T_2$  as above and let us
denote by $x_1,y_1,x_2,y_2$ the elements of $P(\tilde{C},\sigma )$
such that $T_1=\Xi +x_1\cup \Xi +y_1,T_2=\Xi +x_2\cup \Xi +y_2$. There
 are two possible ways of representing the set $\{ x_1,y_1,x_2,y_2\} $
 as union $A\cup B$, where $\# A=\# B=2$, and $A,B$ have one point of
$\{ x_1,y_1\} $ and one point of $\{ x_2,y_2\} $. Namely as :
\begin{equation}\label{e40.1}
\begin{array}{ccl}\{ x_1,x_2\} &\cup &\{ y_1,y_2\} ,\\
\{ x_1,y_2\} &\cup &\{ y_1,x_2\}
  \end{array}
\end{equation}
Taking the sums of the sets in (\ref{e40.1}), using
Corollary~(\ref{c39.1}) and taking into account that $\lambda
_1+\lambda
_2=\mu $ (Section~(\ref{s3})) we see that the extended Prym data
determines the following $P_2$-orbit of quadruples of points in
$P(\tilde{C},\sigma )$, each quadruple being split into a union of
two pairs:
\begin{equation}\label{e405.1}
\begin{array}{c}\{ \{ 2p'_1+\pi ^*p_2-\pi _1^*(2\zeta _1)-\mu
+\rho ,2p''_1+\pi ^*p_2-\pi _1^*(2\zeta _1)-\mu +\rho \} \\
\cup \{ 2p'_2+\pi ^*p_1-\pi _1^*(2\zeta _1)-\mu
+\rho ,2p''_2+\pi ^*p_1-\pi _1^*(2\zeta _1)-\mu +\rho \} \}
\end{array}
\end{equation}
where $\rho \in P_2$. The splitting of the quadruples is consistent
with the action of $P_2$ on $Q$.

Let $f_2^{-1}(f(p_i))=\{ q'_i,q''_i\} $ where $\pi
_2^{-1}(q'_i)=p'_i,\pi _2^{-1}(q''_i)=p''_i,i=1,2$. Let
$f_2^{-1}(f(p_j))=q_j$ for $3\leq j\leq 2g+2$. Since $\pi
_2:C\longrightarrow C_2$ is branched at $\{ q'_1,q''_1,q'_2,q''_2\} $
one has
\begin{equation}\label{e41.2}
2p'_i=\pi _2^*(q'_i)\; ,\; 2p''_i=\pi _2^*(q''_i)
\end{equation}
One has also that
\begin{equation}\label{e41.3}
\pi _1^*(2\zeta _1)=\pi _1^*f^*_1(f(p_i))=2\pi ^*(p_i)
\end{equation}

{\bf Claim} {\it For any $i$ with $1\leq i\leq 2$ and any $j$ with
$3\leq j\leq 2g+2$ there is a point $\rho _{ij}\in P_2$ such that $\mu
=-\pi ^*(p_i-p_j)-\rho _{ij}$.}\\
{\bf Proof.} By Eq.~(\ref{e8.1}) one has $\pi ^*J(C)_2=P_2\cup
(\mu +P_2)$. Since $P(\tilde{C},\sigma )=\pi ^*J(C_2)$ one concludes
by the description of the points of order 2 of the hyperelliptic
Jacobian $J(C_2)$ \cite{mum2} that
$$
P_2=\{ \pi_2^*(S_1-S_2)\mid S_1\cup S_2\subset R_2,S_1\cap
S_2=\emptyset ,\# S_1=\# S_2\}
$$
Using this one easily shows that $\pi ^*(p_i-p_j)\not \in P_2$. q.e.d.

Now, using Eq.~(\ref{e41.2}),(\ref{e41.3}) and the Claim we have
for any $j$ with $3\leq j\leq 2g+2$
\[\begin{array}{cccc}\mbox{}&2p'_1+\pi ^*p_2-\pi _1^*(2\zeta _1)-\mu
+\rho &=&\pi _2^*(q'_1)-\pi ^*(p_2)-\mu +\rho \\
=&\pi _2^*(q'_1)-\pi ^*(p_j)+\rho _{2j}+\rho &=&\pi
_2^*(q'_1-q_j)+\rho _{2j}+\rho
\end{array}\]
We obtain that the $P_2$-orbit $Q$ is equal to
\begin{equation}\label{e425.1}
\begin{array}{c}\{ \{ \pi _2^*(q'_1-q_j)+\rho \; ,\; \pi
_2^*(q''_1-q_j)+\rho\} \\
\cup \{ \pi _2^*(q'_2-q_j)+\rho \; ,\; \pi _2^*(q''_2-q_j)+\rho \} \}
\end{array}
\end{equation}

{\bf Reconstruction of $(C,\eta )$ in the hyperelliptic case, $g\geq 3
.$}\\
    We have a polarized isomorphism $P(\tilde{C},\sigma )\simeq
J(C_2)$. So, by Torelli's theorem \cite{acgh} one reconstructs the
smooth, hyperelliptic curve $C_2$ of genus $g_2=g-1\geq 2$. It has a
unique complete linear system $g^1_2$. Take a point $q\in C_2$ such
that $2q\in g^1_2$. Consider the Abel map $\alpha :C_2\longrightarrow
J(C_2)$ given by $\alpha (x)=cl(x-q)$.
\begin{lem}\label{l435.1}
There is a unique quadruple of $Q$ whose points belong to $\alpha
(C_2)$. Any other quadruple has no points in common with $\alpha
(C_2)$.
\end{lem}
{\bf Proof.} If $q=q_j$ we set $\rho =0$ in (\ref{e425.1}) and see
that
the quadruple
\begin{equation}\label{e435.1}
\{ \{ \alpha (q'_1),\alpha (q''_1)\}\cup \{ \alpha (q'_2),\alpha
(q''_2)\} \}
\end{equation}
is contained in $\alpha (C_2)$. Suppose that for some $\rho \in
J(C_2)_2,\rho \neq0$ one has
$$
q'_1-q_j+\rho \equiv x-q_j.
$$
Multiplying by 2 both sides of this equality we obtain $2q'_1\equiv
2x$. Since $C_2$ has a unique $g^1_2$ and $q'_1\neq x$ for $\rho \neq
0$ we conclude that $2q'_1\in g^1_2$ which is an absurd. This argument
 shows that none of the quadruples of $Q$ different from
(\ref{e435.1})
can have points in common with $\alpha (C_2)$. q.e.d.

Now, we choose a map $f_2:C_2\longrightarrow {\bf P}^1$ of degree 2
and observe that the quadruple of points of $C_2$ defined in the lemma
is transformed by $f_2$ into a set of two points. Furthermore this set
does not depend on the choice of the ramification point $q$ of $f_2$.
Let us denote it by $B_1$. Let $B_2$ be the branch locus of $f_2$.
Then $C$ is isomorphic to the hyperelliptic curve branched at
$B=B_1\cup B_2$ and $\eta \in J(C)_2$ corresponds to this partition of
$B$ \cite{mum2}.

{\bf Reconstruction of $(C,\eta )$ in the case $g=2$.}\\
   Here $P(\tilde{C},\sigma )$ is an elliptic curve $E$. Let $o\in E$
be the zero, let $\varphi  _1,\varphi  _2$ be a basis of $H^0(E,{\cal
O}_{E}(2o))$ and let $f_2=(\varphi  _1:\varphi  _2):E\longrightarrow
{\bf
P}^1$. For any $\rho \in E_2$, if $t_{\rho }:E\longrightarrow E$ is
the translation by $\rho $, there exists $\psi \in PGL(2)$ such that
the following diagram is commutative

\begin{equation}\label{e45.1}
\begin{array}{rlcl}\mbox{}&E&\stackrel{t_{\rho }}{\longrightarrow }
&E\\ f_2&\downarrow &\mbox{}&\downarrow
\hspace{.25cm}f_2\\
\mbox{}&{\bf P}^1&\stackrel{\psi }{\longrightarrow }&{\bf P}^1
\end{array}
\end{equation}
Moreover $\psi $ permutes the branch points of $f_2$. Let $B_2$ be the
branch locus of $f_2$. Take any of the quadruples of $Q$. Each of its
two pairs is invariant under the action of $-id_E$. Thus the image of
the quadruple is a set of two points which we denote by $B_1$. If we
choose another quadruple of $Q$ with image $B'_1$, then
(\ref{e45.1})
shows that there is a $\psi \in PGL(2)$ such that $\psi
(B_1)=B'_1,\psi (B_2)=B_2$. This gives the reconstruction of $(C,\eta
)$, up to isomorphism, as the hyperelliptic curve branched at
$B=B_1\cup B_2$  and $\eta $ as the point of $J(C)_2$ which
corresponds to this partition of $B$.

\section{The bi-elliptic case, $g\geq 4$}\label{s6}
Let $f:C\longrightarrow E$ be a double covering of an elliptic curve
$E$ ramified at $B=\{ x_1,...,x_{2g-2}\} $ and determined by $\delta
\in Pic^{g-1}(E)$ with $\delta ^{\otimes 2}\simeq {\cal O}_{E}(B)$.
Suppose $\eta =f^*(\epsilon )$ for some $\epsilon \in Pic^0(E)_2$.
Then
the unramified covering $\pi :\tilde{C}\longrightarrow C$ determined
by $\eta $ fits into the commutative diagram of Fig.2
\begin{figure}\label{f47.1}
\begin{center}
\begin{picture}(46,54)(0,0)
\put (19,0){\makebox(8,8){$E$}}
\put (0,23){\makebox(8,8){$C_1$}}
\put (19,23){\makebox(8,8){$C$}}
\put (38,23){\makebox(8,8){$C_2$}}
\put (19,46){\makebox(8,8){$\tilde{C}$}}
\put (19,46){\vector(-1,-1){15}}
\put (27,46){\vector(1,-1){15}}
\put (23,46){\vector(0,-1){15}}
\put (23,23){\vector(0,-1){15}}
\put (4,23){\vector(1,-1){15}}
\put (42,23){\vector(-1,-1){15}}
\put (24,15){$f$}
\put (24,38){$\pi $}
\put (6,38){$\pi _1$}
\put (37,38){$\pi _2$}
\put (6,15) {$f_1$}
\put (37,15) {$f_2$}
\end{picture}
\end{center}
\caption{}
\end{figure}
where $deg(f_i)=deg(\pi _i)=2,f_1:C_1\longrightarrow E$ is unramified,
determined by $\epsilon ,f_2:C_2\longrightarrow E$ is ramified at $B$
and is determined by $\delta _2=\delta \otimes \epsilon $. Here we
have the assumptions of Part (iii) of Proposition~(\ref{p26.1}) so the
extended Prym data determines :
\begin{itemize}
\item $E$ as the curve isomorphic to the dual $X\subset {\bf
P}(T_0P)$ of the unique irreducible component of degree $\geq 2$ of
the branch locus $G(R)$ of the Gauss map $G:T^{ns}\longrightarrow
{\bf P}(T_0P)^*$.
\item The points $\{ x_i \mid i=1,...,2g-2\} $ as the duals of the
remaining irreducible components of $G(R)$.
\item $\delta _2\simeq \delta \otimes \epsilon $ as isomorphic to
${\cal O}_{X}(1)$.
\end{itemize}
So, it remains to reconstruct $\epsilon $ which is the
content of the rest of this section.
\begin{lem}\label{l48.1}
Let $T_1=Z_1-\pi ^*\Delta -\mu _1,T_2 = Z_2-\pi ^*\Delta -\mu _2$ be
arbitrary divisors of the orbits $O_1,O_2\subset \mid 2\Xi \mid$,
where $\mu _i\in \lambda _i+P_2$ (Section~(\ref{s3})). Then $T_1,T_2$
are irreducible. Reordering, if necessary, $\{ \lambda _1,\lambda
_2\} $, respectively $\{ O_1,O_2\} ,\{ Z_1,Z_2\} ,\{ T_1,T_2\} $ one
has that the elements $e_1\in T_1,e_2\in T_2$ have the form

(i) $e_1=\pi _1^*(\xi _1)+\pi _2^*(\xi _2)-\pi ^*\Delta-\mu
_1$\\ \noindent where $\xi _1\in C_1,\xi _2\in W_{g-2}(C_2)$ and
$Nm_{f_1}(\xi _1)+Nm_{f_2}(\xi _2)=\delta _2$.

(ii) $e_2=\pi _2^*(\xi _2)-\pi ^*\Delta-\mu _2$\\
where $\xi _2\in W_{g-1}(C_2)$ and $Nm_{f_2}(\xi _2)=\delta
_2$.
\end{lem}
{\bf Proof.} One has to calculate the irreducible components of $Z$
defined in (\ref{e165.1}). One has
$$
H^0(C,K_C\otimes \eta ) \simeq H^0(C,f^*\delta _2)\simeq H^0(E,\delta
_2)\oplus H^0(E,\delta _2\otimes \delta ^{-1})\simeq H^0(E,\delta _2)
$$
Thus $\mid K_C\otimes \eta \mid =f^*\mid \delta _2\mid $. If
$\hat{D}$
is an effective divisor of $\tilde{C}$ such that $Nm_{\pi
}(\hat{D})\in f^*\mid \delta _2\mid $, then
$$
\hat{D}=\pi _1^*E+\pi _2^*F
$$
where $E$ and $F$ are effective divisors of $C_1,C_2$ respectively.
One has $Nm_{\pi }\circ \pi ^*_i=f^*\circ Nm_{f_i}$, so
$$
Nm_{f_1}(E)+Nm_{f_2}(F)\equiv \delta _2
$$
Corollary~(\ref{c19.1}) and a dimension count show that $cl(\hat{D})$
 might be a general element of $Z$ if either $deg(E)=1,deg(F)=g-2$ or
$E=0,deg(F)=g-1$. So, $Z=Z'\cup Z''$ where
$$
Z'=\{ \pi _1^*(\xi _1)+\pi _2^*(\xi _2)\mid \xi _1\in C_1,\xi _2\in
W_{g-2}(C_2),Nm_{f_1}(\xi _1)+Nm_{f_2}(\xi _2)\equiv \delta _2\}
$$
and
$$
Z''=\{ \pi _2^*(\xi _2)\mid \xi _2\in W_{g-1}(C_2),Nm_{f_2}(\xi
_2)\equiv \delta _2\}
$$

   {\bf Claim 1.} {\it $Z'$ is irreducible.}\\
{\bf Proof.} We consider the map $h:C_2^{(g-2)}\longrightarrow E$
defined by $h(D)=\mid \delta _2-Nm_{f_2}(D)\mid $ and the pull-back
diagram

$$
\begin{array}{ccl}X&\longrightarrow &C_1\\
\downarrow &\mbox{}&\downarrow f_1\\
C_2^{(g-2)}&\stackrel{h}{\longrightarrow }&E
\end{array}
$$
Then $Z'$ is the image of $X$ under the map
$$
(D,x)\longmapsto cl(\pi _1^*(x)+\pi _2^*(D))
$$
Now, $X$ might be reducible if $h_*\pi _1(C_2^{(g-2)})$ is contained
in $f_{1*}\pi _1(C_1)$. This is impossible. Indeed,
$f_{2*}:H_1(C_2)\longrightarrow H_1(E)$ is epimorphic since $f_2$ is
ramified. Therefore
$$
(cl\circ f_2^{(g-2)})_*:H_1(C_2^{(g-2)})\longrightarrow
H_1(J_{g-2}(E)) $$
is epimorphic. Composing it with the isomorphism
$J_{g-2}(E)\longrightarrow E$ given by $\xi \mapsto \mid \delta
_2-\xi \mid $ one obtains that $h_*:H_1(C_2^{(g-2)})\longrightarrow
H_1(E)$ is epimorphic. This proves that $X$ and therefore $Z'$ are
irreducible. q.e.d.

   {\bf Claim 2.} {\it $Z''$ is irreducible. }\\
{\bf Proof.} With the same notation as above one considers the
pull-back diagram
$$
\begin{array}{ccl}Y&\longrightarrow &C_2\\
\downarrow &\mbox{}&\downarrow f_2\\
C_2^{(g-2)}&\stackrel{h}{\longrightarrow }&E
\end{array}
$$
Then $Z''$ is the image of $Y$ under the map
$$
(D,y) \longmapsto cl(\pi _2^*(D+y))
$$
In order to prove that $Y$ is irreducible it suffices to verify that
not every component of the branch divisor $h^*(B)$ has even
multiplicity. Now, $h$ can be decomposed as
$$
h=p\circ f_2^{(g-2)}:C_2^{(g-2)}\longrightarrow E^{(g-2)}
\longrightarrow E
$$
where $p$ is the fiber bundle map defined by $p(A)=\mid \delta
_2-A\mid $. Let $x\in B$. Then $\mid \delta _2-x\mid $ is a linear
system of degree $g-2\geq 2$ without base points. Let
$A=p_1+...+p_{g-2}$ be an element with no points in common with $B$.
Let $D\in C_2^{(g-2)}$ with $f_2^{(g-2)}(D)=A$. Then $p^{-1}(x)$ is
smooth at $A$ and $f_2^{(g-2)}$ is nondegenerate at $D$, thus
$h^*(x)$ is a divisor with multiplicity 1. q.e.d.

   Now, $Z$ has two connected components $Z_1$ and $Z_2$, enumerated
as in Proposition~(\ref{p17.1}). So, $Z'\neq Z'', Z_i$ are irreducible
and either $Z_1=Z',Z_2=Z''$ or $Z_1=Z'',Z_2=Z'$. Reordering $\{
Z_1,Z_2\} $ if necessary we can assume that the former case takes
place. q.e.d.
\begin{lem}\label{l51.1}
The singular locus of $Z_1$ has codimension $\geq 2$.
\end{lem}
{\bf Proof.} Consider a divisor of $\tilde{C}$ of the form $\pi
_1^*A+H$ where $H$ is effective, $\pi _1$-simple and $deg(A)\geq 1$.
Then by Corollary~(\ref{c92.2}) one concludes that
\begin{equation}\label{e515.1}
h^0(\tilde{C},\pi _1^*A+H)=h^0(C_1,A)
\end{equation}
Let $\hat{D}=\pi _1^*(x)+\pi _2^*(F)$ where $x\in C_1,F$ is effective
divisor of $C_2$ and $f_1(x)+Nm_{f_2}(F)\equiv \delta _2$.
Proposition~(\ref{p17.1}) and (\ref{e515.1}) show that ${\cal
O}_{\tilde{C}}(\hat{D})$ is a singular point of $Z_1$ if and only if
$F$ is not $\pi _1$-simple. Now, if $\pi _1^*(y)\leq \pi _2^*(F)$,
then one easily checks that $f_2^*(f_1(y))\leq F$. Thus $Sing(Z_1)$
consists of
\begin{equation}\label{e52.1}
cl(\pi _1^*(x+f_1^*(t))+\pi _2^*(G))
\end{equation}
where $x\in C_1,t\in E,G\in C_2^{(g-4)}$ and $Nm_{f_2}(G)\in \mid
\delta _2-f_1(x)-2t\mid $. For any $G\in C_2^{(g-4)}$ there are two
different $\zeta _1\in J_3(C_1)$ such that $Nm_{f_1}(\zeta _1)\equiv
\delta _2-Nm_{f_2}(G)$. Thus the elements of the type (\ref{e52.1})
form a sublocus of $Z_1$ of dimension $\leq g-4$. q.e.d.
\begin{lem}\label{l53.1}
$Sing(Z_2)$ has a unique irreducible component $V$ of codimension 1 in
 $Z_2$. A Zariski open, dense subset of $V$ consists of the elements
\begin{equation}\label{e53.3}
cl(\pi _2^*(f^*_2(x)+G))
\end{equation}
where $x\in E$, and $G$ is an effective, $f_2$-simple divisor of $C_2$
such that $Nm_{f_2}(G)\in \mid \delta _2-2x\mid $.
\end{lem}
{\bf Proof.} We claim that
\begin{equation}\label{e53.1}
dimW^1_{g-1}(C_2)\cap Nm_{f_2}^{-1}(\delta _2)\leq g-4
\end{equation}
This is clear if $C_2$ were not hyperelliptic. If $C_2$ were
hyperelliptic, then $W^1_{g-2}(C_2)=g^1_2+W_{g-3}(C_2)$. This
irreducible variety can not be contained in $Nm_{f_2}^{-1}(\delta
_2)$. Indeed, otherwise its translation would be contained in the
abelian hypersurface $Nm_{f_2}^{-1}(0)$ of $J(C_2)$ which is absurd
since this translation generates $J(C_2)$. By (\ref{e53.1}) we
conclude that the sublocus of $Sing(Z_2)$ :
\begin{equation}\label{e54.1}
\{ \pi _2^*(\xi _2)\mid Nm_{f_2}(\xi _2)=\delta _2,h^0(C,\xi
_2)\geq 2\}
\end{equation}
has codimension $\geq 2$ in $Z_2$.

   Suppose $F$ is an effective, $f_2$-simple divisor of $C_2$ such
that $Nm_{f_2}(F)\in \mid \delta _2\mid $. Assume that $cl(\pi
_2^*F)\in SingZ_2$. By Proposition~(\ref{p17.1}) this is equivalent to
$h^0(\tilde{C},\pi _2^*F)\geq 2$. Then we claim that $h^0(C,F)\geq 2
$, so $\pi _2^*F$ belongs to the locus (\ref{e54.1}). Indeed, since
$\pi
_2:\tilde{C}\longrightarrow C$ is a double unramified covering
corresponding to $f^*_2(\epsilon )\in Pic^0(C_2)_2$ we have
\begin{equation}\label{e53.2}
h^0(\tilde{C},\pi _2^*F)=h^0(C_2,F)+h^0(C_2,f^*_2(\epsilon )(F))
\end{equation}
By Lemma~(\ref{l92.1}) we conclude that $h^0(C_2,f^*_2(\epsilon
)(F))=0$. So, $cl(\pi _2^*(F))$ belongs to $Sing(Z_2)$ if and only if
$cl(F)\in W^1_{g-1}(C_2)$.

    Now, suppose that $F=f_2^*(x)+G$ where $x\in E,G$ is effective and
$Nm_{f_2}(G)\in \mid \delta _2-2x\mid $. Let $t_{\epsilon }(x)$ be the
translation of x by $\epsilon $. Then by Eq.~(\ref{e53.2}) we have
\begin{equation}\label{e54.3}
h^0(\tilde{C},\pi
_2^*F)=h^0(C_2,f_2^*(x)+G)+h^0(C_2,f_2^*(t_{\epsilon }(x))+G)
\end{equation}
Thus $h^0(\tilde{C},\pi ^*F)\geq 2$ and $cl(\pi ^*F)\in Sing Z_2$. The
sublocus of $SingZ_2$
$$
V=\{ cl(\pi _2^*(f_2^*(x)+G))\mid x\in E,G\geq 0,Nm_{f_2}(G)\in \mid
\delta _2-2x\mid \}
$$
is the image of $X$ where $X$ is defined by the pull-back diagram
\begin{equation}\label{e55.1}
\begin{array}{ccl}X&\longrightarrow &E\\
\downarrow &\mbox{}&\downarrow \beta \\
C_2^{(g-3)}&\stackrel{\alpha }{\longrightarrow }&J_2(E)
\end{array}
\end{equation}
and $\alpha (G)=cl(\delta _2-Nm_{f_2}(G)),\beta (x)=cl(2x)$. The same
argument as in Claim 1 of Lemma~(\ref{l48.1}) shows that $X$ is
irreducible. This implies that $V$ is irreducible as well.
Corollary~(\ref{c92.2}) and Eq.~(\ref{e54.3}) show that for $F=f
_2^*(x)+G$ one has $h^0(\pi _2^*F)=2$ if and only if $G$ is
$f_2$-simple. Thus the points (\ref{e53.3}) form a Zariski open, dense
subset of $V$.

   Finally, $dimX=g-3$ and we claim that the map $X\longrightarrow V$
given by $(G,x)\longmapsto cl(\pi _2^*(f_2^*(x)+G))$ is of degree 2,
hence $dimV=g-3$. Indeed, let $\sigma _2:\tilde{C}\longrightarrow
\tilde{C}$ be the involution which interchanges the sheets of $\pi
_2$. Then for any $\pi _2^*(\xi _2)\in V$ with $h^0(\tilde{C},\pi
_2^*(\xi _2))=2$ according to Eq.~(\ref{e54.3}) there are exactly two
$\sigma _2$-invariant divisors in $\mid \pi _2^*(\xi _2)\mid $ namely
\begin{equation}\label{e56.1}
\begin{array}{cc}
\pi _2^*(f_2^*(x)+G)\; ,&\pi _2^*(f_2^*(t_{\epsilon }(x))+G)
\end{array}
\end{equation}
with $x,G$ as in the lemma. q.e.d.

   Let
$$
S_2=\{ F\in C_2^{(g-2)}\mid Nm_{f_2}(F)\in \mid \delta _2\mid \}
$$
One has a surjective map
$$
\varphi  =cl\circ \pi _2^*:S_2\longrightarrow Z_2
$$
{}From the proof of Claim 2 of Lemma~(\ref{l48.1}) we see that $S_2$
is irreducible. Moreover $deg\varphi  =1$ since $h^0(\tilde{C},L)=1$
for
any sufficiently general $L\in Z_2$. Let us consider the Stein
factorization \cite{hart}
$$
\varphi  =\psi \circ \alpha :S_2\longrightarrow \Gamma \longrightarrow
Z_2
$$
where $\psi $ is a finite map and $\alpha $ has connected fibers. Let
$$
W_1=\{ F\in S_2\mid h^0(C_2,F)\geq 2\}\; ,\;
W_2=\{ f_2^*A+E\in S_2\mid A\geq 0,E\geq 0,deg(A)\geq 2\}
$$
One has $codim_{S_2}(W_1)\geq 1$ and $codim_{S_2}(W_2)\geq 2$. Let
$S_2^0=S_2\backslash (W_1\cup W_2)$ and let $\Gamma ^0=\alpha(S^0_2)$.
\begin{lem}\label{l58.1}
The points of $S^0_2,\Gamma^0$ are nonsingular in $S_2,\Gamma $
respectively, $codim_{\Gamma }(\Gamma \backslash \Gamma ^0)\geq 2$ and
the map $\alpha:S^0_2\longrightarrow \Gamma ^0$ is an isomorphism. Let
$n:N\longrightarrow Z_2$ be the normalization of $Z_2$. Then there
exists a finite map $\beta :N\longrightarrow \Gamma $ such that
$n=\psi \circ \beta $. If $N^0=\beta ^{-1}(\Gamma ^0)$ then
$codim_N(N\backslash N^0)\geq 2$ and $\beta :N^0\longrightarrow
\Gamma ^0$ is an isomorphism
\end{lem}
{\bf Proof.} The points of $S^0_2$ are nonsingular by
Proposition~(\ref{p91.1}). For any $x\in S^0_2$ one has $\# \varphi
^{-1}(\varphi  (x))=2$ as we have shown in the proof of
Lemma~(\ref{l53.1}).
Thus the map $\alpha :S^0_2\longrightarrow \Gamma ^0$ is bijective.
The map $\varphi  $ is nondegenerate at the points of $S^0_2$. Indeed
$\varphi  = cl\circ \pi _2^*=\pi _2^*\circ cl$, the map
$cl:S_2\longrightarrow J_{g-1}(C_2)$ is nondegenerate at any $F\in
S^0_2$ since $h^0(F)=1$, and the map $\pi
_2^*:J_{g-1}(C_2)\longrightarrow J_{2g-2}(\tilde{C})$ is obviously
nondegenerate. We conclude that $\alpha :S^0_2\longrightarrow \Gamma
^0$ is an isomorphism. One has $codim_{\Gamma }(\Gamma \backslash
\Gamma ^0)\geq 2$ since $codim_{Z_2}\varphi  (W_1)\geq 2$ by
(\ref{e53.1}).
The rest of the lemma is clear by the universal property of the
normalization. q.e.d.

{\bf Reconstruction of $(C,\eta )$ in the bi-elliptic case, $g\geq
4$.}

\noindent
In the beginning of this section we have seen how to reconstruct up to
isomorphism $E$ and the covering $f_2:C_2\longrightarrow E$. Let
$T_i\in O_i\subset  \mid 2\Xi \mid ,i=1,2$. We have proved above that
$T_i$ are irreducible and just one of $T_i$ has a singular locus of
codimension 1. Reordering, if necessary, as in Lemma~(\ref{l48.1}) we
can assume that this divisor is $T_2$. We can identify $T_2$ and $Z_2$
 by translation. Let $n:N\longrightarrow T_2$ be the normalization of
$T_2$. Let $R=n^{-1}(V)$. The Zariski open subset $R^0=R\cap N^0$
is dense in $R$ since $codim_N(N\backslash N^0)\geq 2.$ By the
irreducibility of $X$ in (\ref{e55.1}) one concludes that
$\alpha^{-1}\circ \beta (R^0)$ and $R$ are irreducible as well.  We
have an isomorphism $f^*:\mid \delta _2\mid \longrightarrow \mid
K_C\otimes \eta \mid $. Consider the Gauss map
$G:Z_2^{ns}\longrightarrow \mid K_C\otimes \eta\mid $. Then for every
$F\in S_2$ with $h^0(\tilde{C},\pi _2^*F)=1$ one has
$$
G(\varphi
  (F))=f^*(Nm_{f_2}(F))
$$
Shrinking $N_0$ from Lemma~(\ref{l58.1}) we can assume that the
following properties hold
\begin{itemize}
\item $codim_N(N\backslash N_0)\geq 2$.
\item The composition $C\circ n$ can be extended to a regular map on
$N_0$.
\item Every point of $\alpha ^{-1}\circ \beta (R^0)$ has the form
$f_2^*(x)+A$ where $x$ is not a branch point of $f_2$, $A$ is reduced
and $f_2$-simple, and $x\not \in Supp(A).$
\end{itemize}
    Now, we define a rational map
$$
\gamma : R\longrightarrow E
$$
as follows. For every $L\in R_0$ the hyperplane $G\circ n(L)$ belongs
to the unique irreducible component of degree $>1$ of the branch locus
of the Gauss map $G:T\longrightarrow {\bf P}(T_0P)^*$, namely $E^*$.
By the conditions above this hyperplane is tangent to a unique point
of $E$. We denote this point by $\gamma (L)$. Now, let $L=\beta
^{-1}\circ \alpha (f_2^*(x)+A)\in R_0$. Then $n^{-1}(n(L))=\{ L,L'\} $
where
\begin{equation}\label{e62.1}
L'=\beta ^{-1}\circ \alpha (f_2^*(t_{\epsilon  }(x))+A)
\end{equation}
according to (\ref{e56.1}). This shows that the map
$n:R\longrightarrow
V\subset Sing(T_2)$ is of degree 2. By (\ref{e62.1}) the corresponding
involution $\tau ^*:{\bf C}(R)\longrightarrow {\bf C}(R)$ of the field
of rational functions on $R$ transforms $\gamma ^*{\bf C}(E)$ into
itself and $\tau ^*:\gamma ^*{\bf C}(E)\longrightarrow \gamma ^*{\bf
C}(E)$ is induced by the translation map $t_{\epsilon
}:E\longrightarrow E$. This gives the reconstruction of $\epsilon \in
J(E)_2$ and completes the reconstruction of $(C,\eta )$ from the
extended Prym data.

\section{The case $g=3$}\label{s7}
Let $a\in \tilde{C}$. We define the Abel-Prym map $\phi
:\tilde{C}\longrightarrow P(\tilde{C},\sigma )$ by
$$
\phi (x)=cl(x-a-\sigma (x-a))
$$

\begin{lem}\label{l635.1}
Suppose $g\geq 2$. The following alternative takes place
\begin{list}{(\roman{bean})}{\usecounter{bean}}
\item $\phi $ maps $\tilde{C}$ isomorphically onto its image $\phi
(\tilde{C})$.
\item The map $\phi :\tilde{C}\longrightarrow \phi (\tilde{C})$ has
degree 2.
\end{list}
The second case occurs if and only if

{\rm (*)} $C$ is hyperelliptic and $\eta \simeq {\cal O}_{C}(p_1-p_2)$
for some $p_1,p_2\in C$.\\
\noindent Here $\phi (\tilde{C})\simeq C_2$ and $\phi =\pi _2$ (see
Fig.1) via this isomorphism.
\end{lem}
{\bf Proof.} Suppose $\phi (x)=\phi (y)$ for some $x\neq y$. Then
$x+\sigma y\equiv y+\sigma x$, thus $\tilde{C}$ is hyperelliptic. It
has a unique $g^1_2$, so $\sigma (g^1_2)=g^1_2$. Let $\sigma _1$ be
the hyperelliptic involution of $\tilde{C}$. Then $\sigma \neq \sigma
_1$ and we claim that $\sigma $ and $\sigma _1$ commute. Indeed, for
any $z\in \tilde{C}$
\[ \begin{array}{cc}
\sigma z+\sigma _1(\sigma z)\in g^1_2\; ,&\sigma (z+\sigma _1z)\in
g^1_2
\end{array} \]
Thus $\sigma \sigma
_1z=\sigma _1\sigma z$. Let $\sigma _2=\sigma \sigma _1$. Let
$C_1=\tilde{C}/\sigma _1$ and let $\overline{\sigma
}:C_1\longrightarrow C_1$ be the involution induced by $\sigma $. Then
$\overline{\sigma }$ has two fixed points since $C_1\simeq {\bf P}^1$.
Thus $\pi :\tilde{C}\longrightarrow C$ fits into the commutative
diagram of Fig.1 and condition (*) holds.
    If $\phi $ were degenerate at some point $x\in \tilde{C}$, then
$\pi (x)$ would be a base point of $\mid K_C\otimes \eta\mid $, thus
condition (*) holds according to Lemma~(\ref{l21.1}).
   Conversely, suppose condition (*) holds. Then by the argument above
$\phi (x)=\phi (y)$ and $x\neq y$ if and only if $y$ belongs to the
divisor $\sigma (x+\sigma _1x)$. Thus $y=\sigma _2(x)$. q.e.d.

   Further we suppose that $g=3$. Let $T_i\in O_i\subset \mid 2\Xi
\mid ,i=1,2.$ The divisors $T_i$ are reduced, connected curves
according to Proposition~(\ref{p17.1}). Let $S=Nm^{-1}(\mid
K_C\otimes \eta \mid )\subset \tilde{C}^{(4)} $ and let
$Z=Nm^{-1}(K_C\otimes \eta )\cap W_4(\tilde{C})\subset
J_4(\tilde{C})$. Both $S$ and $Z$ break into two disjoint, connected
components $S=S_1\cup S_2,Z=Z_1\cup Z_2$. We enumerate so that
$cl(S_i)=Z_i$ and $T_i$ is translation of
$Z_i,i=1,2$.
\begin{lem}\label{l66.1}
The curves $T_1,T_2$ are both singular if and only if condition (*)
of Lemma~(\ref{l635.1}) holds. If only one of $T_i$ is singular then
the
nonsingular one is a translation of $\phi (\tilde{C})$.
\end{lem}
{\bf Proof.} If condition (*) holds, then both $T_1$ and $T_2$ are
reducible and hence singular by Corollary~(\ref{c39.1}). Since $T_i$
is
a translation of $Z_i,i=1,2$ we can work with $Z_i$. Suppose that
condition (*) does not hold and $Z'\in \{ Z_1,Z_2\} $ is a singular
curve with a singular point $L$. Let $\{ Z',Z''\} =\{ Z_1,Z_2\} $.
Then
$$
X=\{ L+x-\sigma x\mid x\in \tilde{C}\} \subset Z''
$$
Clearly $X$ is a translation of $\phi (\tilde{C})$. According to
Lemma~(\ref{l635.1}) $X$ is isomorphic to $\tilde{C}$ and $X$ is
algebraically equivalent to $2\Xi $ \cite{mas}. Thus $X=Z''$. q.e.d.

{\bf Reconstruction of $(C,\eta ) $ in the case $g=3$.}

{\bf Case 1.} {\it Both $T_1,T_2$ are singular.}\\
According to Lemma~(\ref{l66.1}) we are in the situation of
Section~(\ref{s5}) where a procedure for the reconstruction of
$(C,\eta )$ was described.

{\bf Case 2.} {\it Just one of the curves $T'\in \{ T_1,T_2\} $ is
singular.}\\
We take the other curve $T''$. It is isomorphic to
$\tilde{C}$ according to Lemmas~(\ref{l66.1}) and (\ref{l635.1}). The
involution
$-id_P:T''\longrightarrow T''$ coincides with $\sigma $ via this
isomorphism. We thus reconstruct $(\tilde{C},\sigma )$.

{\bf Case 3.} {\it $T_1$ and $T_2$ are nonsingular.}\\
Consider the involutions $\sigma _i:T_i\longrightarrow T_i$
induced by $-id_P$, let $C_i=T_i/\sigma _i$ and let $\pi
_i:T_i\longrightarrow C_i$ be the factor maps, $i=1,2$. By
Proposition~(\ref{p17.1}) the map $cl:S_i\longrightarrow Z_i$ is
bijective. Thus $S_i$ are nonsingular since $Z_i$ are nonsingular.
Using Proposition~(\ref{p91.1}) one checks that the nonsingularity of
$S_i$ implies that $\sigma ^{(4)}:S_i\longrightarrow S_i$ is without
fixed points, thus $\sigma _i:T_i\longrightarrow T_i$ is without fixed
points as well, $i=1,2$. Consider the Gauss maps
$G_i:T_i\longrightarrow {\bf P}(T_0P)^*={\bf P}^1$. By
Eq.~(\ref{e26.1})
one shows that $G_i=f_i\circ \pi _i$ where $f_i:C_i\longrightarrow
{\bf P}^1$ are maps of degree 4 and concludes that the maps
$G_i:T_i\longrightarrow {\bf P}^1$ are obtained from $\varphi _{
K_C\otimes \eta }\circ \pi :\tilde{C}\longrightarrow {\bf P}^1$ by the
tetragonal construction of Donagi \cite{don},\cite{don1}.
   Now, take the pair $(T_1,\sigma _1)  $ and apply the tetragonal
construction to $G_1:T_1\longrightarrow {\bf P}^1$ (ibid.). One
obtains two 8-sheeted coverings $g_i:X_i\longrightarrow {\bf P}^1$
with involutions $\tau _i:X_i\longrightarrow X_i$ such that $g_i\circ
\tau _i=g_i$. For one of them, e.g. $X_2$, there is an isomorphism of
the coverings
$$
\begin{array}{rlcl}\mbox{}&X_2&\stackrel{\psi _2}{\longrightarrow
}&T_2\\ g_2&\downarrow &\mbox{}&\downarrow
\hspace{.25cm}G_2\\ \mbox{}&{\bf P}^1&=&{\bf P}^1
\end{array}
$$
such that $\psi _2\circ \tau _2=\sigma _2\circ \psi _2$. Then the
remaining pair $(X_1,\tau _1)$ is isomorphic to $(\tilde{C},\sigma )$.

\noindent {\sc Institute of Mathematics, Bulgarian Academy of
Sciences,}\\
\noindent Acad. G. Bonchev Str. bl. 8, {\sc Sofia 1113 Bulgaria}

\end{document}